\documentclass[12pt,twoside]{amsart}
\usepackage{graphicx}
\usepackage{amsmath}
\usepackage{amssymb}
\usepackage{amscd}
\usepackage[francais]{babel}
\usepackage[latin1]{inputenc}
\usepackage{amsmath,amssymb,amsthm}
\usepackage[french]{minitoc}
\usepackage[all]{xy}
\usepackage{pstricks}
\usepackage{graphicx}
\usepackage{enumerate,paralist}
\graphicspath{./}
\usepackage{latexsym,amscd,amsfonts}

\theoremstyle{plain}
\title{Bornes pour la r\'egularit\'e de Castelnuovo-Mumford des sch\'emas non lisses}
\author{Amadou Lamine Fall}
\address{D\'epartement de Math\'ematiques, Universit\'e Cheikh Anta Diop, Dakar, S\'en\'egal et Institut Math\'ematique de Jussieu, paris VI}

\email{fall@math.jussieu.fr}
\date{}
\newtheorem{theoreme}{Théorème}[section]
\newtheorem{lemme}[theoreme]{Lemme}
\newtheorem{proposition}[theoreme]{Proposition}
\newtheorem{corollaire}[theoreme]{Corollaire}

\theoremstyle{definition}
\newtheorem{definition}[theoreme]{Définition} 

\newcommand{\Z}{{\mathcal Z}}
\newcommand {\pp}{\mathbb{P}}

\newcommand {\oo}{{\mathcal O}}
\newcommand {\yy}{{\mathcal Y}}

\newcommand{\fm}{{\mathfrak m}}

\newcommand{\fp}{{\mathfrak p}}
\newcommand{\fb}{{\mathfrak b}}

\newcommand{\HH}{H_{\mathfrak m}}

\newcommand{\lra}{\longrightarrow}

\newcommand{\codim}{\operatorname{codim}}

\newcommand{\reg}{\operatorname{reg}}

\newcommand{\TR}{\operatorname{Tor}^R}

\newcommand{\indeg}{\operatorname{indeg}}
\newcommand{\proj}{\operatorname{Proj}}

\newcommand{\sing}{\operatorname{Sing}}

\begin{document}

\pagestyle{plain}
\begin{abstract}
Nous montrons dans cet article  des bornes pour la r\'egularit\'e de Castelnuovo-Mumford d'un  sch\'ema admettant des singularit\'es, en fonction des degr\'es des \'equations d\'efinissant le sch\'ema, de sa dimension et de la dimension de son lieu singulier. Dans le cas o\`u les singularit\'es sont isol\'ees, nous améliorons la borne fournie par Chardin et Ulrich dans \cite{CU} et dans le cas g\'en\'eral, nous \'etablissons une borne doublement exponentielle en la
dimension du lieu singulier.
\end{abstract}

\maketitle

\section{Introduction}
Nous \'etudions  des bornes pour la r\'egularit\'e de Castelnuovo-Mumford d'un  sch\'ema, admettant des singularit\'es en fonction des degr\'es des \'equations d\'efinissant le sch\'ema, de sa dimension et de celle de son lieu singulier. Ce travail est une continuation des travaux de Bertram-Ein-Lazarsfeld d'une part et Chardin-Ulrich d'autre part.

Soit $R=k[X_{0},\ldots X_{n}]$, un anneau de polyn\^omes sur le corps $k$, $I\subset R$ un id\'eal homog\`ene engendr\'e par des \'el\'ements homog\`enes de degr\'es au plus $D$. Soit $X=\proj(R/I)$  le sch\'ema projectif sur  $k$ d\'efini par $I$, $d$ sa dimension et $r>0$ sa codimension. Soit $\delta$ la dimension du lieu singulier de $X$ (avec la convention $\delta=-1$ si $X$ est lisse). 

Si la caract\'eristique du corps $k$ est nulle, $X$ purement de codimension $r$ et $\delta\leq 0$, Bertram-Ein-Lazarsfeld \cite{BEL} dans le cas lisse et Chardin-Ulrich \cite{CU} dans le cas o\`u les singularit\'es sont isol\'ees, ont montr\'e la borne suivante: $$\reg(I_{X})\leq r(D-1)+1.$$
Dans cette note nous \'etablissons les bornes suivantes, qui sont valables en toute caract\'eristique:\medskip
\begin{equation}
 \reg(I_{X})\leq (\dim X)!(r(D-1)-1)+1\quad\quad  si \quad \delta\leq 1 
\end{equation}
\begin{equation}
\reg(I_{X})\leq \lambda D^{(n-\delta)2^{\delta-2}}\quad\quad si \quad \delta\geq 2 
\end{equation}\medskip

o\`u $\lambda$ est explicit\'e dans le th\'eor\`eme 4.3 et ne d\'epend que de $n$, $d$ et $\delta$.\medskip

Pour \'etablir nos bornes, nous proc\'edons en deux \'etapes. Dans la premi\`ere \'etape, nous \'etablissons des bornes pour la r\'egularit\'e des  sch\'emas  dont les singularit\'es sont rationnelles et localement intersection compl\`etes, hors d'un nombre fini de points. Nous utilisons pour cela la m\'ethode de Chardin et Ulrich. Cette m\'ethode  d\'ecrite dans \cite{CU}, utilise des techniques de liaison, une r\'ecurrence sur la dimension et une version am\'elior\'ee du th\'eor\`eme d'annulation de Kodaira.
Dans la deuxi\`eme \'etape, on se ram\`ene au cas \'etudi\'e dans la premi\`ere \'etape, en utilisant un th\'eor\`eme de Bertini et une r\'ecurrence introduite par Caviglia et Sbarra dans \cite{CS} et d\'evelopp\'ee par Chardin, Fall et Nagel dans \cite{CFN}.\vskip 1cm
{\bf Remerciements.} {\it Je remercie Marc Chardin, qui m'a propos\'e le sujet et  aid\'e par des discussions et remarques pertinentes \`a r\'ealiser ce texte.}\bigskip

\section{R\'esultats pr\'eliminaires}
Dans cette section nous rappelons les r\'esultats et d\'efinitions qui seront utilis\'es dans les autres parties du texte.

Soit $R=k[X_{0},\ldots,X_{n}]$ un anneau de polyn\^omes sur un corps $k$, $\fm=(X_{0},\ldots,X_{n})$ et $M$ un $R$-module gradu\'e de type fini. 
Posons $b_{i}^{R}(M)=\max\{ \mu / \TR_{i}(M,k)_{\mu}\neq 0 \}$ si $\TR_{i}(M,k)\neq 0$ et $b_{i}^{R}(M)=-\infty $ sinon. Notons $\HH^{i}(M)$
le $i$-i\`eme module de cohomologie locale de $M$ \`a support dans $\fm$, $a_{i}^{R}(M)=\max\{ \mu /
\HH^{i}(M)_{\mu}\neq 0 \}$ si $\HH^{i}(M)\neq 0$ et $a_{i}^{R}(M)=
-\infty $ sinon. On rappelle que l'on d\'efinit la regularit\'e de Castelnuovo-Mumford de $M$ par: 
$$
\reg(M) = \max_{i} \{a_{i}^{R}(M)+i\}
= \max_{i} \{b_{i}^{R}(M)-i\} .
$$

Le nombre minimal de g\'en\'erateurs d'un module $M$ sur un anneau local (ou d'un $R$-module gradu\'e $M$) est not\'e $\mu (M)$.

\subsection{Bornes pour la régularité des sch\'emas en dimension au plus un.}

Soit  $I=(f_{1},\ldots,f_{s})\subset R$ un id\'eal homog\`ene de hauteur $n-1$,  engendr\'e par $s$ formes de degr\'es $d_{1}\geq d_{2}\geq \ldots \geq d_{s}$ et soit  $X=\proj(R/I)$ le sch\'ema projectif d\'efini par $I$ et $C$ la composante de dimension 1 de $X$. On a le r\'esultat suivant qui am\'eliore \cite[Proposition 2.1 et Proposition 2.2]{CU}:

\begin{theoreme}\cite[Proposition 5.12 et Th\'eor\`eme 5.13]{Ch1}\\
Avec les notations et les hypoth\`eses ci-dessus on a:
$$\quad\quad \reg(R/I)= \sum_{i=1}^{s}(d_{i}-1)\quad si \quad s=n-1.$$   Si $C$ est r\'eduite,
  $$\quad \reg(R/I)\leq 2(\sum_{i=1}^{n-1}(d_{i}-1)-1)+d_{n} \quad si \quad s\geq n .$$
Si de plus,  $\mu(I_{\fp})\leq n$ pour tout id\'eal premier $\fp \supset I_{C}$ tel que $dim(R/\fp)=1$, alors $$\reg(R/I)\leq \sum_{i=1}^{\min\{s,n+1\}}(d_{i}-1).$$
\end{theoreme}\medskip

\subsection{Sur les singularit\'es et sur un th\'eor\`eme de Bertini.}
Soit $X\subset \pp_{k}^{n}$ un sch\'ema projectif sur un corps $k$ et $x\in X$ un point ferm\'e de $X$.
\begin{definition}
On dit que $x$ est un point non singulier de $X$ si l'anneau local $\oo_{X,x}$ est r\'egulier. Si $\oo_{X,x}$ n'est pas r\'egulier  on dit que $x$ est un point singulier de $X$.
\end{definition}\medskip
Le lieu singulier du sch\'ema $X$, $\sing X$, est l'ensemble de ses points singuliers. Notons que $\sing(X)$ est un ferm\'e de $X$.
\begin{definition}
Soit $A$ un anneau noeth\'erien. On dit que $A$ v\'erifie la condition $R_{\ell}$, si $A_{\fp}$ est r\'egulier pour tout id\'eal premier $\fp$ de $A$ de hauteur au plus $\ell$.
\end{definition}

\begin{definition}
On dit qu'un sch\'ema $X$ v\'erifie la condition $R_{\ell}$ si pour tout $x\in X$, l'anneau $\oo_{X,x}$ v\'erifie la condition $R_{\ell}.$
\end{definition}\medskip
 Le r\'esultat suivant donne une version du th\'eor\`eme de Bertini qui est valable en toute caract\'eristique.
\begin{theoreme}\cite[3.4.14]{FOW} \\
Soit $X\subset \pp_{k}^{n}$ un sch\'ema  projectif sur un corps infini $k$. Si $X$ est r\'egulier (respectivement normal, r\'eduit, v\'erifie $R_{\ell}$), alors pour tout hyperplan g\'en\'eral $H$, $X \cap H$ est r\'egulier (respectivement normal, r\'eduit, v\'erifie $R_{\ell}$).
\end{theoreme}\medskip

\begin{corollaire}
Soit $X\subset \pp_{k}^{n}$ un sch\'ema projectif et $H$ un hyperplan g\'en\'eral. Si $\dim(\sing(X))=s$, alors $\dim(\sing(X\cap H))=s-1$.
\end{corollaire} 
Dans ce qui suit, nous rappelons les d\'efinitions de singularit\'es de type rationnel et singularit\'es $F$-rationelles.
\begin{definition}
Soit $X$ un sch\'ema de type fini sur un corps $k$, on dit que  $X'\buildrel{\pi}\over{\lra} X$ est une d\'esingularisation de $X$ si $\pi$ est un morphisme propre birationnel et si $X'$ est lisse sur $k$.
\end{definition}\medskip
\begin{definition}\cite{K}\\
Soit $X$ un sch\'ema essentiellement de type fini sur un corps $k$ de caract\'eristique z\'ero, $X'\buildrel{\pi}\over{\lra} X$ une d\'esingularisation de $X$, soit $R^{i}\pi_{*} \oo_{X'}$ les images directes sup\'erieures du faisceau $\oo_{X'}$. On dit que $X$ poss\`ede  des singularit\'es rationnelles si $X$ est normal et si $R^{i}\pi_{*} \oo_{X'}=0$ pour tout $i>0$.
\end{definition}\medskip

\begin{definition}
Un anneau $R$ de caract\'eristique premi\`ere est $F$-rationnel si tout id\'eal de $R$ engendr\'e par un syst\`eme de param\`etres est \'etroitement clos (tightly closed). Un sch\'ema est $F$-rationnel si tous ses anneaux locaux sont $F$-rationnel.
\end{definition}\medskip

La notion de $F$-rationalit\'e s'etend aux anneaux essentiellement de type fini sur un corps $k$ de caract\'eristique z\'ero  (voir \cite[D\'efinition 4.1]{S}). Elle coïncide avec la notion de singularit\'e rationnelle.\\
Nous adoptons la terminologie de Chardin et Ulrich dans \cite{CU}:

\begin{definition}
  On dit qu'un anneau $R$ est de type rationnel, s'il  est de caract\'eristique $p$ et  $F$-rationnel ou bien s'il  est essentiellement de type fini sur un corps de caract\'eristique z\'ero et \`a singularit\'es rationnelles. Un sch\'ema $X$ est de type rationnel si ses anneaux locaux sont de type rationnel.
\end{definition}\medskip

\section{Bornes pour la r\'egularit\'e des sch\'emas  en dimension au moins  $2$}

Soit  $R$ un anneau de polyn\^omes sur un corps, le th\'eor\`eme suivant  am\'eliore  le th\'eor\`eme \cite[4.7]{CU}  de Chardin et Ulrich:\medskip
\begin{theoreme}
Soit $R=k[X_{0},\ldots,X_{n}]$ un anneau de polyn\^omes sur un corps $k$, $I=(f_{1},\ldots,f_{s})\subset R$ un id\'eal non nul engendr\'e par des formes de degr\'es $d_{1}\geq \ldots \geq d_{s}\geq 2$,  $X=\proj(R/I)$. Posons $r=\codim X$ et  $\sigma=\sum_{i=1}^{r}(d_{i}-1)$.
 On suppose qu'il existe un sch\'ema $\Z \subset X$ de dimension z\'ero, tel que pour tout $x\in X-\Z$, $X$ est localement intersection compl\`ete en $x$ et $\oo_{X,x}$ est de type rationnel.

 Alors
 $\reg(R/I)\leq \sigma$ si $R/I$ est de Cohen-Macaulay, et sinon
$$
\begin{array}{lcll}
(0)&\quad\quad &\reg(R/I)\leq \sum_{i=1}^{n+1}(d_{i}-1) & \hbox{si} \ \dim X \leq 0, \\
(1)&\quad\quad &\reg(R/I)\leq 2(\sigma-1)+d_{n} & \hbox{si} \ \dim X=1,\\
(2) &\ & \reg(R/I)\leq (\dim X+1)!(\sigma-1) & \hbox{si} \ \dim X\geq 2.\\
\end{array}
$$
\end{theoreme} \bigskip

\begin{proof}
Si $R/I$ est de Cohen-Macaulay, $I$ est contenu dans une intersection compl\`ete $\fb$ de degr\'es $d_{1}\geq d_{2}\geq \ldots \geq d_{r}$. D'apr\`es \cite[4.1(a)]{CU}, 
$$
\reg(R/I)=\sigma-\indeg \left(\frac{\fb :I}{\fb}\right) \leq\sigma.
$$

 Le $(0)$ d\'ecoule de \cite[3.3]{Ch1} (voir aussi \cite{Sj}) et le $(1)$ d\'ecoule du th\'eor\`eme 2.1.\\
  
Dans ce qui suit nous utilisons la m\'ethode et les notations de  \cite[4.4 et 4.7]{CU}.

Soit $r=\codim X$, si $r=1$, on  peut supposer que $n\geq 3$, il s'agit  de montrer que $\reg(R/I)\leq n!(d_{1}-1)$. Comme les  $f_{i}$ ont un facteur commun $h$, on  posons $f_{i}=hf'_{i}$ et $I'=(f'_{1},\ldots,f'_{r})$, l'id\'eal $I'$ est de codimension  $r'\geq 2$ et $\reg(I)=\deg(h)+\reg(I')$.

Si $r'\geq n$, on a d'apr\`es $(0)$,  $\reg(R/I')\leq (n+1)(d_{1}-\delta-1)$, o\`u $\delta=\deg(h)$, donc 
\begin{eqnarray*}
\reg(R/I) &\leq & (n+1)(d_{1}-\delta-1)+\delta \\
& \leq & (n+1)(d_{1}-1)\\
& \leq & n!(d_{1}-1).
\end{eqnarray*}

Si $r'\leq n$, le fait que $X$ soit localement intersection compl\`ete hors d'un sch\'ema $\Z$ de dimension, implique que $\dim(R/I')\leq 2$. Donc, les in\'egalit\'es $(1)$ et $(2)$ appliqu\'ees \`a $I'$ donne la borne pour $I$.

Dans tout ce qui suit on suppose $r\geq 2$. Posons $$a_{i,j}=\sum_{\arrowvert \mu \arrowvert=d_{i}-d_{j}}U_{i,j,\mu}X^{\mu}\quad 1\leq i\leq r \quad et  \quad r+1\leq j\leq s, $$ o\`u les $U_{i,j,\mu}$ sont des variables, $X^{\mu}= X_{0}^{\mu_{0}}\ldots X_{n}^{\mu_{n}}$ et $\arrowvert \mu \arrowvert = \mu_{0}+\ldots +\mu_{n}$. Consid\'erons la matrice $A=(a_{i,j})$ et d\'efinissons  $\alpha_{1},\ldots,\alpha_{r}$ par: $$(\alpha_{1},\ldots,\alpha_{r})= \begin{pmatrix}
I_{r,r} & A 
\end{pmatrix}
\begin{pmatrix}
f_{1}\\.
\\. \\. \\ f_{s}
\end{pmatrix},$$ o\`u $I_{r,r}$ est la matrice identit\'e d'ordre $r$. Posons $K=k(U_{i,j,\mu})$, $R'=R\otimes_{k}K$,   $J=(\alpha_{1},\ldots,\alpha_{r})R':IR'$, $\yy=\proj(R'/IR'+J).$ 
Pour $c=\dim R-1$ et $c'=\dim R-2$,  le th\'eor\`eme \cite[4.4(d)]{CU} montre que $\yy$ est un sch\'ema de type rationnel et localement intersection compl\`ete hors d'un sch\'ema de dimension z\'ero. D'apr\`es \cite[1.7(iii) ]{CU}, $\yy$ coincide avec $\yy'=\proj(R'/IR'+(J)_{\leq \sigma})$ hors d'un nombre fini de points. De plus, comme $\yy'$ est localement intersection compl\`ete hors d'un nombre fini de points, il existe  $d=\dim X$ formes $\beta_{1},\ldots,\beta_{d} \in J$ de degr\'es au plus $\sigma$ telles que $\yy''=\proj(R'/(I,\beta_{1},\ldots,\beta_{d}))$ coïncide avec $\yy'$ hors d'un nombre fini de points. Ainsi $\yy''$ coïncide avec $\yy$ hors d'un nombre fini de points. En posant $J''=(\alpha_{1},\ldots,\alpha_{r},\beta_{1},\ldots,\beta_{d})$, on a la suite exacte suivante:
$$0\lra R'/IR'\cap J''\lra R'/IR'\oplus R'/J''\lra R'/IR'+J''\lra 0.$$
De cette suite on d\'eduit que 

\begin{eqnarray*}  (*)\quad \reg(R/I) & =& \reg(R'/IR')  \\ & \leq & \max\{\reg(R'/IR'\cap J''),\reg(R'/IR'+J'')\}.
\end{eqnarray*} 

Comme $IR'\cap J''=(\alpha_{1},\ldots,\alpha_{r})$ est une intersection compl\`ete de codimension $r$, on a $$\reg(R'/IR'\cap J')=\sigma.$$

Puisque  $\yy''$ coïncide avec $\yy$ hors d'un nombre fini de points, $\yy''$ est un sch\'ema de dimension $d-1$, de type rationnel et localement intersection compl\`ete hors d'un sch\'ema de dimension z\'ero. Nous allons en d\'eduire $(2)$ par r\'ecurrence sur la dimension $d$ de $X$.\\ 

Pour $d=2$,  $\yy''$ est d\'efini par des \'equations de degr\'es  $\sigma\geq \sigma \geq d_{1}\geq \ldots \geq d_{r}$.  Comme $\dim \yy''=1$, on a 
\begin{eqnarray*}
\reg(R/I) & \leq & \reg(R'/IR'+J'')\\
& =& 2(d(\sigma-1)+\sum_{i=1}^{r-d}((d_{i}-1)-1)+d_{r-d+1}\\
& =& 2(2(\sigma-1)+\sum_{i=1}^{r-2}(d_{i}-1)-1)+d_{r-1}\\
& \leq & 2(2(\sigma-1)+(\sigma-1))\\
& = & 6(\sigma-1).
\end{eqnarray*}\medskip
Pour  $d\geq 3$, $\yy''$ est d\'efini par des \'equations de degr\'es $\underbrace{\sigma \geq \ldots \geq \sigma}_{d \quad fois }\geq  d_{1}\geq \ldots \geq d_{r-d}$. D'apr\`es l'hypoth\`ese de r\'ecurrence, si $r\geq d$, on a  
\begin{eqnarray*}
\reg(R'/IR'+J'') & \leq & d!(d(\sigma-1)+\sum_{i=1}^{r-d}(d_{i}-1)-1)\\
& \leq & d!(d(\sigma-1)+\sigma-1)\\
& =& (d+1)!(\sigma-1).
\end{eqnarray*}

 et si $r\leq d$, on a
\begin{eqnarray*}
\reg(R'/IR'+J'') & \leq & d!(r(\sigma-1))\\
& \leq & dd!(\sigma-1)\\
& = & (d+1)!(\sigma-1).
\end{eqnarray*}

On en d\'eduit dans tous les cas que, $$ \reg(R/I)\leq (d+1)!(\sigma-1).$$
\end{proof}\bigskip
\begin{corollaire}
Soit $X\subset \pp^{n}$ un sch\'ema projectif sur un corps $k$, de codimension $r>0$, d\'efini par des \'equations de degr\'es au plus $D\geq 2$. On suppose qu'il existe un sch\'ema $\Z \subset X$ de dimension z\'ero, tel que $\forall x\in X-\Z$, $X$ est localement intersection compl\`ete en $x$ et $\oo_{X,x}$ est de type rationnel. Alors,\bigskip
$$
\begin{array}{lcll}

(0)&\quad\quad &\reg(R/I)\leq (n+1)(D-1) & \hbox{si} \ \dim X \leq 0, \\
(1)&\quad\quad &\reg(R/I)\leq (2r+1)(D-1)-1  & \hbox{si} \ \dim X=1,\\
(2) &\ & \reg(R/I)\leq (\dim X+1)!(r(D-1)-1) & \hbox{si} \ \dim X\geq 2.\\
\end{array}
$$
\end{corollaire}\bigskip

\section{Bornes pour la r\'egularit\'e des sch\'emas  non lisses}\bigskip
Soit $R$ un anneau de polyn\^omes sur un corps, $I$ un id\'eal  gradu\'e de $R$, engendr\'e par des \'el\'ements de degr\'es au plus $D$ et soit $l_{1},\ldots,l_{s+1}$ des formes lin\'eaires g\'en\'erales. Posons $M=R/I$,  $M_{i}:=M/(l_{1},\ldots ,l_{i})M$, $i = 0,\ldots,s$ et $K_{i}:=\ker (M_{i}\buildrel{\times l_{i+1}}\over{\lra}M_{i}[1])$. En appliquant les r\'esultats de la preuve de \cite[3.2]{CFN} avec $R=S$ et $J=(0)$, on a 
\begin{eqnarray*}
r_{i} &=& \max \{ b_{1}^{R}(M)-2,\ b_{0}^{R}(M)+\max \{ 1, b_{0}^{S}(J)\} -2,\ \reg(M_{i})\}\\
& = & \max \{b_{1}^{R}(M)-2,\  b_{0}^{R}(M)-1,\ \reg(M_{i})\}\\
& = & \max \{ D-2,\reg(M_{i})\}
\end{eqnarray*}

 On a ainsi le lemme suivant comme cas particulier de \cite[3.2]{CFN}:\bigskip

\begin{lemme} Avec les notations ci-dessus, on suppose que les 
  $K_{i}$ sont de longueur finie pour tout $i$.
En posant $Q_{i}=\max\{\reg(M_{i}),\lambda(K_{i}), D-2 \}+1$ pour $0\leq i\leq s$, on a  $$Q_{i}\leq Q_{i+1}^{2} \quad \forall i = 0,\ldots,s-1.$$  En particulier,
$$
\reg (M)\leq Q_{s}^{2^{s}}.
$$ 
\end{lemme} \bigskip

 Le lemme 4.1   et la proposition suivante nous permettent d'obtenir des bornes pour la r\'egularit\'e du sch\'ema $X$.\medskip

\begin{proposition}
Soit $R=k[X_{0},\ldots,X_{n}]$ un anneau de polyn\^omes sur un corps $k$, avec $n\geq 2$, $\fm=(X_{0},\ldots,X_{n})$, $I\subset R$ un id\'eal engendr\'e par des formes de degr\'es au plus $D$ et $X= \proj(R/I)$. On d\'esigne par $d$ la dimension de $X$ et par $\delta$ celle de son lieu singulier.
Soit $l_{1},\ldots,l_{s}$,  avec $s\geq \delta$ des formes lin\'eaires g\'en\'erales, on pose  $S=R/I$, $S_{i}:=S/(l_{1},\ldots ,l_{i})S$, $i=0,\ldots,s$, $K_{i}:=\ker (S_{i}\buildrel{\times l_{i+1}}\over{\lra}S_{i}[1])$ et $X_{\delta}=\proj(R/I+(l_{1},\ldots,l_{\delta}))$. Supposons qu'il existe un sch\'ema $\Z_{\delta} \subset X_{\delta}$ de dimension z\'ero, tel que $\forall x\in X_{\delta}-\Z_{\delta}$, $X_{\delta}$ est localement intersection compl\`ete en $x$ et $\oo_{X_{\delta},x}$ est de type rationnel.
 Alors,

$$\lambda(K_{\delta-1})\leq D^{r}\left( \frac{(r(D-1)(d+1)!)^{d+1-\delta}}{(d+1-\delta)!}\right)$$
\end{proposition}\bigskip
\begin{proof}
La suite exacte  $$0\lra K_{\delta-1}\lra S_{\delta-1}\buildrel{\times l_{\delta}}\over{\lra}S_{\delta-1}[1]\lra S_{\delta}\lra 0$$

 donne une suite exacte en cohomologie locale 
$$0\lra K_{\delta-1}\lra H^{0}_{\fm}(S_{\delta-1})\buildrel{\times l_{\delta}}\over{\lra} H^{0}_{\fm}(S_{\delta-1})[1]\lra H^{0}_{\fm}(S_{\delta})\lra\ldots,$$ qui montre que  
\begin{eqnarray*}
\lambda(K_{\delta-1}) & \leq & \lambda(H^{0}_{\fm}(S_{\delta-1}))-\lambda(H^{0}_{\fm}(S_{\delta-1}))+\lambda(H^{0}_{\fm}(S_{\delta})).\\
& = & \lambda(H^{0}_{\fm}(S_{\delta-1}))
\end{eqnarray*}
Ainsi,  nous avons 
\begin{eqnarray*} \lambda(K_{\delta-1}) & \leq & \sum_{\nu=0}^{\reg(S_{\delta})}\lambda(H^{0}_{\fm}(S_{\delta}))_{\nu}\\
& \leq & \sum_{\nu=0}^{\reg(S_{\delta})}H_{S_{\delta}}(\nu).
\end{eqnarray*} 

L'id\'eal $I+(l_{1},\ldots,l_{\delta})$ contient un id\'eal $J_{\delta}$ engendr\'e par une suite r\'eguli\`ere de degr\'es $\underbrace{D,\ldots,D}_{r\quad fois},\underbrace{1,\ldots,1}_{\delta \quad fois}$. Donc 
\begin{eqnarray*}
H_{S_{\delta}}(\nu)& \leq & H_{R/J_{\delta}}(\nu)\\
& = & \sum_{i_{1}=0}^{D-1}\ldots \sum_{i_{r}=0}^{D-1}{{\nu+d-\delta-(i_{1}+\ldots + i_{r})}\choose{d-\delta}}\\
& \leq & {{\nu+d-\delta}\choose{d-\delta}}D^{r},
\end{eqnarray*}
ainsi,

\begin{eqnarray*}
\lambda(K_{\delta-1})& \leq & D^{r}\sum_{\nu=0}^{\reg(S_{\delta})}{{\nu+d-\delta}\choose{d-\delta}}\\
& =& D^{r}{{\reg(S_{\delta})+d+1-\delta}\choose{d+1-\delta}}
\end{eqnarray*} 
Comme $\reg(S_{\delta})\leq (d+1)!(r(D-1)-1)$, on a 

\begin{eqnarray*}
\lambda(K_{\delta-1})& \leq & D^{r}{{(d+1)!(r(D-1)-1)+d+1-\delta}\choose{d+1-\delta}}\\
& \leq & D^{r}{{r(d+1)!(D-1)}\choose{d+1-\delta}}\\
& \leq & D^{r}\left( \frac{(r(D-1)(d+1)!)^{d+1-\delta}}{(d+1-\delta)!}\right)
\end{eqnarray*}
\end{proof}\bigskip
\begin{theoreme}
Soit $R=k[X_{0},\ldots X_{n}]$ un anneau de polyn\^omes sur le corps $k$, $I\subset R$ un id\'eal homog\`ene, $X=\proj(R/I)$, $d=\dim X$    et  $\delta$ la dimension du lieu singulier de $X$.  On suppose que $I$ est  engendr\'e par des \'el\'ements de degr\'es au plus $D\geq 2$.\\ Alors on a
$$
\reg(R/I)\leq C(n,d,\delta)D^{(n+1-\delta)2^{\delta-1}}, $$ \vskip 0.5cm
 o\`u
$C(n,d,\delta)=\left(\frac{((n-d)(d+1)!)^{d+1-\delta}}{(d+1-\delta)!}\right)^{2^{\delta-1}}$.

\end{theoreme}\bigskip
\begin{proof}
Soit   $l_{1},\ldots,l_{\delta}$ des formes lin\'eaires g\'en\'erales, on pose $$S=R/I,\quad S_{i}:=S/(l_{1},\ldots ,l_{i})S,\quad i=0,\ldots,\delta \quad  et \quad K_{i}:=\ker (S_{i}\buildrel{\times l_{i+1}}\over{\lra}S_{i}[1]).$$  Si $H_{1},\ldots,H_{\delta}$ d\'esignent les hyperplans d\'efinis par les $l_{i}$, le sch\'ema $\Z = X\cap H_{1}\cap \ldots \cap H_{\delta}$ v\'erifie les conditions  du corollaire 3.2. Comme $(d+1-\delta)! < (d+1)!$, on a $$\reg(S_{\delta})\leq (d+1)!(r(D-1)-1).$$  D'apr\`es la preuve de \cite[3.2]{CFN}, nous  avons  $$\reg(S_{\delta-1})\leq \reg(S_{\delta})+\lambda(K_{\delta-1}).$$ 

 Donc, 
\begin{eqnarray*}
 Q_{\delta-1} & = & \max\{\reg(S_{\delta-1}),\lambda(K_{\delta-1}), D-2 \}+1\\
& \leq & \max \{\reg(S_{\delta})+\lambda(K_{\delta-1}), D-2\}+1\\
& \leq & (d+1)!(r(D-1)-1)+\left( \frac{(r(D-1)(d+1)!)^{d+1-\delta}}{(d+1-\delta)!}\right)D^{r}\\
& \leq &  \left(\frac{(rD(d+1)!)^{d+1-\delta}}{(d+1-\delta)!}\right)D^{r}\\
& \leq & \left(\frac{(r(d+1)!)^{d+1-\delta}}{(d+1-\delta)!}\right)D^{r+d+1-\delta}\\
& = & \left(\frac{((n-d)(d+1)!)^{d+1-\delta}}{(d+1-\delta)!}\right)D^{n+1-\delta}.
\end{eqnarray*}
Ainsi d'apr\`es lemme 4.1, nous avons 
\begin{eqnarray*}
\reg(R/I) & \leq & Q_{\delta-1}^{2^{\delta-1}}\\
& \leq & \left(\frac{((n-d)(d+1)!)^{d+1-\delta}}{(d+1-\delta)!}\right)^{2^{\delta-1}}D^{(n+1-\delta)2^{\delta-1}}.
\end{eqnarray*}

\end{proof}\bigskip
En utilisant les th\'eor\`emes 3.1 et 4.3, les r\'esultats de Bertram-Ein-Lazarsfeld \cite{BEL} et les r\'esultats de Chardin-Ulrich \cite{CU}, nous avons le th\'eor\`eme suivant:\medskip
\begin{theoreme}
Soit $X$  un sch\'ema projectif sur un corps $k$, de dimension $d$   et  de codimension $r>0$. Soit $\delta$ la dimension du lieu singulier de $X$ et $I_{X}$ l'id\'eal satur\'e d\'efinissant $X$. On suppose que $X$ est d\'efini par des \'equations de degr\'es au plus $D\geq 2$.\vskip 0.5cm
 $1)$ Si $\delta =-1$ ou si $\delta =0$ et la caract\'eristique de $k$ est nulle, alors 
$$ \reg(I_{X})\leq r(D-1)+1.$$ \vskip 0.5cm
$2)$ Si $\delta\leq 1$, alors $\reg(I_{X})\leq (\dim X)!(r(D-1)-1)+1.$\vskip 0.5cm
$3)$ Si $\delta\geq 2$ alors, 
$$\reg(I_{X})\leq C(n,d-1,\delta-1)D^{(n-\delta)2^{\delta-2}},$$
o\`u $C(n,d,\delta)$ est  d\'efinie dans le th\'eor\`eme 4.3.
\end{theoreme}\medskip
\begin{proof}
La conclusion du $1)$ d\'ecoule des th\'eor\`emes \cite[ 4 (i)]{BEL} et de \cite[ 0.1]{CU}. 

Soit $R=k[X_{0},\ldots,X_{n}]$ et $I_{X}=I^{sat}$ l'id\'eal satur\'e d\'efinissant le sch\'ema $X$. Soit $l$  une forme lin\'eaire g\'en\'erale et $H$ l'hyperplan d\'efini par $l$, d'apr\`es \cite[2.5]{CM}, on a 
\begin{eqnarray*}
(*)\quad \quad \quad \reg(R/I_{X}) &=& \max\{ \reg(R/I_{X\cap H}), a_{1}(R/I_{X})+1\}\\
& \leq & \max\{ \reg(R/I+(l)), a_{1}(R/I_{X})+1\}.
\end{eqnarray*}
D'autre part la suite exacte, avec $\lambda (K)$ fini,  $$ 0\lra K \lra (R/I)(-1)\buildrel{\times l}\over{\lra} R/I \lra R/I+(l) \lra 0,$$ donne une suite exacte en cohomologie locale $$ \ldots \lra H^{0}_{\fm}(R/I+(l))\lra H^{1}_{\fm}(R/I)(-1)\lra H^{1}_{\fm}(R/I)\lra H^{1}_{\fm}(R/I+(l))\lra \dots $$ qui montre que $H^{1}_{\fm}(R/I)_{\mu}\simeq H^{1}_{\fm}(R/I)_{\mu-1}=0$ pour $\mu \geq \reg(R/I+(l))$. On en d\'eduit que $$\reg(R/I+(l))\geq a_{1}(R/I)+1= a_{1}(R/I_{X})+1,$$ et par suite d'apr\`es $(*)$ on a 
\begin{eqnarray}
\reg(R/I_{X}) \leq \reg(R/I+(l)).
\end{eqnarray}\medskip

Si $\delta\leq 1$ le sch\'ema, $\Z=\proj(R/I+(l))$ est \`a singularit\'es isol\'ees. D'apr\`es le $(2)$ du th\'eor\`eme 3.1 on a,  
\begin{eqnarray*}\reg(R/I+(l)) & \leq &  (\dim \Z +1)!(r(D-1)-1)\\ 
& = & (\dim X)!(r(D-1)-1).
\end{eqnarray*}

On en d\'eduit que 
\begin{eqnarray*}\reg(I_{X}) &=& \reg(R/I_{X})+1\\
& \leq & \reg(R/I+(l))+1\\
 & \leq &  (\dim \Z +1)!(r(D-1)-1)+1\\ 
& = & (\dim X)!(r(D-1)-1)+1.
\end{eqnarray*}

Si $\delta\geq 2$, l' in\'egalit\'e $(3)$ ci-dessus et le th\'eor\`eme 4.3 appliqu\'e \`a $R/I+(l)$ donne la borne annonc\'ee.

\end{proof}

\end{document}